\documentclass[11pt]{article}
\usepackage{amsmath}
\usepackage{amssymb}
\usepackage{amsfonts}
\usepackage{mathrsfs}
\usepackage{mathtools}
\usepackage{enumerate}
\usepackage[usenames]{color}
\usepackage[numbers]{natbib}
\usepackage{enumitem} 
\usepackage{amscd} 

\usepackage[normalem]{ulem} 
\usepackage{comment}
\usepackage{bbm}

\usepackage[pdftex]{graphicx}
\usepackage[all]{xy}

\usepackage{hyperref}
\hypersetup{colorlinks=true,linkcolor=blue,
citecolor=blue} 

\usepackage{amsthm}

\usepackage{latexsym}
\newsavebox{\toy}
\savebox{\toy}{\framebox[0.65em]{\rule{0cm}{1ex}}}
\newcommand{\QED}{\usebox{\toy}\end{demo}}

\numberwithin{equation}{section}

\newtheorem{thm}{Theorem}[section]
\newtheorem{lem}[thm]{Lemma}
\newtheorem{prop}[thm]{Proposition}

\theoremstyle{definition}
\newtheorem{definition}{Definition}[section]
\newtheorem{rem}[definition]{Remark}

\newcommand{\bd}{\begin{displaymath}}
\newcommand{\ed}{\end{displaymath}}
\newcommand{\N}{{\mathbb{N}}}
\newcommand{\Z}{{\mathbb{Z}}}

\newcommand{\R}{{\mathbb{R}}}
\newcommand{{\rd}}{\R^d}

\newcommand{\lan}{\langle}
\newcommand{\ran}{\rangle}

\renewcommand{\b}{\beta}

\newcommand{\del}{\delta}

\newcommand{\e}{\varepsilon}
\newcommand{\ve}{\epsilon}

\newcommand{\kF}{\mathcal{F}}

\newcommand{\kZ}{{\mathcal Z}}

\newcommand{\dd}{\text{\rm d}}             

\newcommand{\dis}{\displaystyle}

\allowdisplaybreaks[4]

\newcommand{\MN}[1]{\textcolor{red}{#1}}

\author{
   Makoto Nakashima\thanks{Research supported by JSPS KAKENHI Grant Numbers JP18H01123, JP18K13423, 	JP22K03351, JP23K22399.}
}

\begin{document}
\title{A note on the asymptotics of the free energy of $1+1$ dimensional directed polymers in random environment at high temperature}

\date{}
\pagestyle{myheadings}
\markboth{}{}

\maketitle


\begin{abstract}
The author gave the sharp asymptotic behavior of the free energy of $1+1$ dimensional directed polymers in random environment(DPRE) as the  inverse temperature $\b\to 0$ under the assumption that random environment satisfies a certain concentration inequality in \cite{Nak19}, \begin{align*}
\lim_{\b\to0}\frac{1}{\beta^4}F(\beta)=-\frac{1}{6}.
\end{align*}
In this paper, we obtain the same asymptotics without using the concentration inequality.
\end{abstract}

\vspace{1em}




\vspace{1em}{\bf Key words:}  Directed polymers,  Free energy, Universality, Continuum directed polymer.

\vspace{1em}
We denote by  $(\Omega, {\cal F},P )$ a probability space. We denote by $P[X]$ the expectation of random variable $X$ with respect to $P$.  Let $\mathbb{N}_0=\{0,1,2,\cdots\}$, $\mathbb{N}=\{1,2,3,\cdots\}$, and $\mathbb{Z}=\{0,\pm 1,\pm 2,\cdots\}$.  
Let $C_{x_1,\cdots,x_p}$ or $C(x_1,\cdots,x_p)$ be a non-random constant which depends only on the parameters $x_1,\cdots,x_p$. 



 \section{Introduction and main result}

In this paper, we consider the directed polymers in random environment(DPRE) on $\Z^d$. DPRE is defined by a Gibbs measure $\mu_N$ on the path space of simple random walk on $\Z^d$, which is described in terms of simple random walk $(S,P_S^0)$ and random potential $\{\eta(n,x):n\in \N, x\in \Z^d\}$. More precisely, $\mu_N$ is given by \begin{align}
\mu_N(dS)=\frac{1}{Z_N}\exp\left(\beta\sum_{k=1}^N \eta (k,S_k)\right)P^0_S(dS),\label{eq:Gibbsmeas}
\end{align}
where $Z_N=Z_{\beta,N}:=P_S^0\left[\exp\left(\beta\sum_{i=1}^N \eta(i,S_i)\right)\right]$ is so-called partition function, which makes $\mu_N$ be a probability measure.

Under $\mu_N$, the shape of polymer, which is described by the realization of  $S$, is interacted by the potential. In particular, it is known that the phase transition of the shape occurs: If the inverse temperature $\beta\geq 0$ is sufficiently small and $d\geq 3$, then the shape of polymer is diffusive as $N\to\infty$\cite{Bol89,Sin95,AZ96,SZ96,CSY04,Var06,Jun23,Lac24,Jun24}. On the other hand, if $\beta\geq 0$ is large enough or $d=1,2$, then the polymer is localized\cite{CSY03,CH02,CH06,JL24}. This phase transition is characterized in terms of the free energy $F(\beta)$. 


The reader may refer to the monograph \cite{Com17} and \cite{Zyg24} for the further results.

\subsection{Model and main result}

To define \eqref{eq:Gibbsmeas} precisely, we introduce some random variables.
\begin{itemize}
\item (Random environment) Let $\{\eta(n,x):(n,x)\in\mathbb{N}\times \mathbb{Z}^d\}$ be $\mathbb{R}$-valued i.i.d.\,random variables defined on a probability space $(\Omega_Q,\kF_Q,Q)$ be a probability space with $\lambda(\beta)=\log Q[\exp\left(\beta \eta(n,x)\right)]\in\mathbb{R}$ for any $\beta\in\mathbb{R}$, and \begin{align}
Q[\eta(n,x)]=0,\quad Q[\eta(n,x)^2]=1\label{eq:normalization},
\end{align} 
where  \eqref{eq:normalization} is a convenient normalization. 
\item (Simple random walk) Let $S$ be a simple random walk on $\mathbb{Z}$ starting from $x=(x_1,\dots,x_d)\in\mathbb{Z}^d$ defined on a probability space $(\Omega_S,\kF_S,P_S^x)$ with $x_1+\dots+x_d\in 2\mathbb{Z}$. We write $P_S=P_S^0$ for simplicity.
\end{itemize}

We also assume that $S$ and $\eta$ are independent under $P_S^x\otimes Q$.


Then, it is clear that \begin{align*}
Q\left[Z_{\beta,N}(\eta)\right]=\exp\left(N\lambda(\beta)\right)
\end{align*}
for any $\beta\in\mathbb{R}$.

The normalized partition function is defined by \begin{align}
 W_{\beta,N}(\eta)&=\frac{Z_{\beta,N}(\eta)}{Q\left[Z_{\beta,N}(\eta)\right]}=P_S\left[\prod_{k=1}^N\zeta_{k,S_k}(\beta,\eta)\right],\label{partfn}
 \end{align}
 where we write for each $(n,x)\in\mathbb{N}\times \mathbb{Z}$\begin{align*}
 \zeta_{n,x}(\beta,\eta)=\exp\left(\beta \eta(n,x)-\lambda(\beta)\right).
 \end{align*}
 
Then, the following limit exists $Q$-a.s.\,and $L^1(Q)$ \cite{CSY03,CY06}:
\begin{align}
F(\beta)&=\lim_{N\to \infty}\frac{1}{N}\log W_{\beta,N}(\eta)=\lim_{N\to\infty}\frac{1}{N}Q\left[\log W_{\beta,N}(\eta)\right]=\sup_{N\geq 1}\frac{1}{N}Q\left[\log W_{\beta,N}(\eta)\right].\label{free}
\end{align}
The limit $F(\beta)$ is a non-random constant called the quenched free energy. Jensen's inequality implies that \begin{align*}
F(\beta)\leq \lim_{N\to\infty}\frac{1}{N}\log Q\left[W_{\beta,N}(\eta)\right]=0.
\end{align*}
It is known that $F(\beta)<0$ if $\beta\not=0$ when $d=1,2$ \cite{CV06,Lac10} and $F(\beta)=0$ for sufficiently small $|\beta|$ when $d\geq 3$. The author obtained  the asymptotics of $F(\beta)$ in the  high temperature regime ($\beta\to 0$) under the technical condition:\begin{align}
&F(\beta)\sim -\frac{\beta^{4}}{6},\ \ \text{if }d=1\label{eq:asymd1}
\intertext{\cite{Lac10, Wat12, AY15,Nak19} and}
&\log |F(\beta)|\sim -\frac{\pi}{\beta^2},\ \ \text{if }d=2\notag
\end{align}
\cite{Lac10,BL17,Nak14}.

We remark that \eqref{eq:asymd1} is closely related to the KPZ equation. Indeed,  the author proved in \cite{Nak19} that under the concentration assumption \eqref{A} below   \begin{align*}
\lim_{\beta\to 0}\frac{1}{\beta^4}F(\beta)=F_\kZ(\sqrt{2})=\lim_{T\to\infty}\frac{1}{T}P_\kZ\left[\log \int_{\mathbb{R}}\mathcal{Z}_{\sqrt{2}}^x(T,y)dy\right]=-\frac{1}{6},
\end{align*}
where $\mathcal{Z}^x_{\beta}(t,y)$ is the unique mild solution to the stochastic heat equation \begin{align*}
\partial \mathcal{Z}=\frac{1}{2}\Delta \mathcal{Z}+\beta\mathcal{Z}\dot{\mathcal{W}},
\end{align*}
with  the initial condition $\dis \lim_{t\to 0 }\mathcal{Z}(t,y)dy=\delta_x(dy)$ and $\mathcal{W}$ is a time-space white noise and $P_\kZ$ is the law of $\kZ_\beta^x$.
We write \begin{align*}
\kZ_{\beta}^x(t)=\int_\mathbb{R}\kZ_{\beta}^x(t,y)dy
\end{align*}
and $\kZ_\beta(t)=\kZ_\beta^0(t)$ for simplicity.


Our main result removes the concentration assumption \eqref{A}.

\begin{thm}\label{thm:main}
Suppose $d=1$ and \eqref{eq:normalization} is satisfied. Then, we have \begin{align*}
\lim_{\beta\to 0}\frac{1}{\beta^4}F(\beta)=F_\kZ(\sqrt{2})=\lim_{T\to\infty}\frac{1}{T}P_\kZ\left[\log \int_{\mathbb{R}}\mathcal{Z}_{\sqrt{2}}^x(T,y)dx\right]=-\frac{1}{6},
\end{align*}
\end{thm}

\begin{rem}\label{rem:1}
In \cite{Nak19}, $(\eta,Q)$ is supposed to satisfy the following concentration inequality:

{There exist $\gamma\geq 1$, $C_1,C_2\in(0,\infty)$ such that for any $n\in\mathbb{N}$ and  for any convex and $1$-Lipschitz function $f:\mathbb{R}^n\to \mathbb{R}$,}
\begin{align}Q\left(|f(\omega_1,\cdots,\omega_n)-Q[f(\omega_1,\cdots,\omega_n)]|\geq t\right)\leq C_1\exp\left(-C_2t^{\gamma}\right),\label{A}
\end{align}
where $1$-Lipschitz means $|f(x)-f(y)|\leq |x-y|$ for any $x,y\in\mathbb{R}^n$ and $\omega_1,\cdots,\omega_n$ are i.i.d.~random variables with the marginal law $Q(\eta(n,x)\in dy)$.

In \cite{Nak19}, the concentration assumption \eqref{A} was used to prove the inequality \begin{align}
\varliminf_{\beta\to 0}\frac{1}{\beta^4}F(\beta)\geq \lim_{T\to\infty}\frac{1}{T}P_\kZ\left[\log \int_{\mathbb{R}}\mathcal{Z}_{\sqrt{2}}^x(T,y)dx\right].\label{eq:lower}
\end{align}
and the opposite inequality was proved without the concentration assumption. Thus, the main task in this paper is to prove \eqref{eq:lower} without \eqref{A}.

\end{rem}

\begin{rem}
Our proof uses the reflection principle  and Chernoff bounds for simple random walk. We expect that they are not essential. Actually, the reflection principle is used to allows us to bound the difference of partition functions simply,  and  Chernoff bounds is used to verifies that the partition function with underlying random walk going far does not contribute to the free energy.  
\end{rem}

\begin{rem}
In \cite{CTT17}, the universal asymptotics of free energy are also discussed for the disordered pinning model with $\alpha\in (\frac{1}{2},1)$ under the assumption of concentration inequality of disorder(Remark \ref{rem:1}). The proof uses the  concentration inequality to obtain the upper bound and the lower bound of the free energy.

We find that the arguments in this paper and \cite{Nak19} can be applied to disordered pinning model, and hence we can remove the concentration inequality assumption in \cite{CTT17} (\cite{Nak25+}). 
\end{rem}

\subsection{Outline of the proof}
As mentioned in Remark \ref{rem:1}, it is enough to show \eqref{eq:lower}. 

Here, we give an outline of the proof of \eqref{eq:lower}. 

First, we introduce the coarse grained space-time lattice $\mathbb{L}_{CG}$  with size $\beta^{-2}\times T\beta^{-4}$ (Subsection \ref{2.1}). Then, we assign the open-closed states to the time-space bonds due to the configuration of $\eta$ in each bond $\lan (I,X),(I+1,Y) \ran$ (Subsection \ref{2.2}). Roughly, we assign an ``open" state to a bond when a ``partition function  $W_{\beta,T}((I,X),(I+1,Y))$" in the bond is bounded from below by $e^{-\left(\frac{1}{6}+o(1)\right)T}$. If the probability that each space-time bond is open is sufficiently close to $1$, then we can find an infinite oriented open path with positive probability so that the partition function $W_{\beta,TN\beta^{-4}}\geq e^{-TN\left(\frac{1}{6}+o(1)\right)}$ for any $T,N\geq 1$ with positive probability. This implies that $\frac{1}{\beta^4}F(\beta)\geq -\frac{1}{6}+o(1)$. 

 An idea of coarse grained oriented percolation was uesd in \cite{AY15}, where most parts are devoted to the asymptotics of free energy of directed polymers with defect line. (The reader may refer  to \cite{KKN21} for the direct discussion for DPRE on the strongly recurrent graph.)

\section{Proof of Theorem \ref{thm:main}}\label{2}

\subsection{Coarse grained argument deduced from DPRE}\label{2.1}

For fixed $T\in [1,\infty)$ and $\delta>0$, we consider the segments \begin{align*}
B^{(T,n,\delta)}_j=[(j\delta T-1)\sqrt{n},(j \delta T+1)\sqrt{n}],\quad j\in\Z,
\end{align*}  
which plays a role of sites in the coarse graining lattice defined below.

We  introduce the   coarse grained space-time lattice:
\begin{align*}
\mathbb{L}_{CG}&=\{(I,X)\in \mathbb{N}_0\times \Z:0\leq |X|\leq I, I-|X|\in 2\mathbb{N}_0 \}
\end{align*}
and we identify $ITn\times B_{X}^{(T,n,\delta)}$ with $(I,X)\in \mathbb{L}_{CG}$. Also, we assign a bond $\langle (I,X),(J,Y)\rangle$ for $(I,X),(J,Y)\in \mathbb{L}_{CG}$ if and only if \begin{align*}
|I-J|=1\text{ and }|X-Y|=1.
\end{align*} 

To bound the partition function, we restrict the trajectories of the simple random walk $S$.

Then, we equip a tube $T^{(T,n,\delta)}(I,X)$ around a vertex $(I,X)$ which contains an edge $\langle (I,X),(I+1,Y)\rangle$ in $\N_0\times \Z$;
\begin{align*}
&\mathtt{T}^{(T,n,\delta,L)}_{(I,X)}=\left\{(k,l)\in \N_0\times \Z: \begin{array}{l}
\dis ITn< k\leq (I+1)Tn,\\
\dis  \left|l
-X\delta T\sqrt{n}\right|\leq LT\sqrt{n}
\end{array}\right\},
\end{align*} 
and for a bond $\langle (I,X),(I+1,Y)\rangle$, we define a set $\Omega_{(I,X,Y)}^{(T,n,\delta)}$  of trajectories  of random walk $S$ 
from $ITn$ to $(I+1)Tn$ by 
\begin{align*}
&\Omega_{(I,X,Y)}^{(T,n,\delta,L)}=\left\{s=\{(i,s_i)\}_{i=ITn}^{(I+1)Tn}\subset \Omega_S: s_{ITn}=x\in B_X^{(T,n,\delta)},s_{(I+1)Tn}\in B_{Y}^{(T,n,\delta)}, s\subset \mathtt{T}_{(I,X)}^{(T,n,\delta,L)}\right\},
\end{align*}
where $L>2\delta+1$ is chosen sufficiently large later. 
We will omit the superscripts $T,n,\delta,L$ when it is clear from the context.

\begin{figure}[htbp]
\begin{center}
\includegraphics[width=3in]{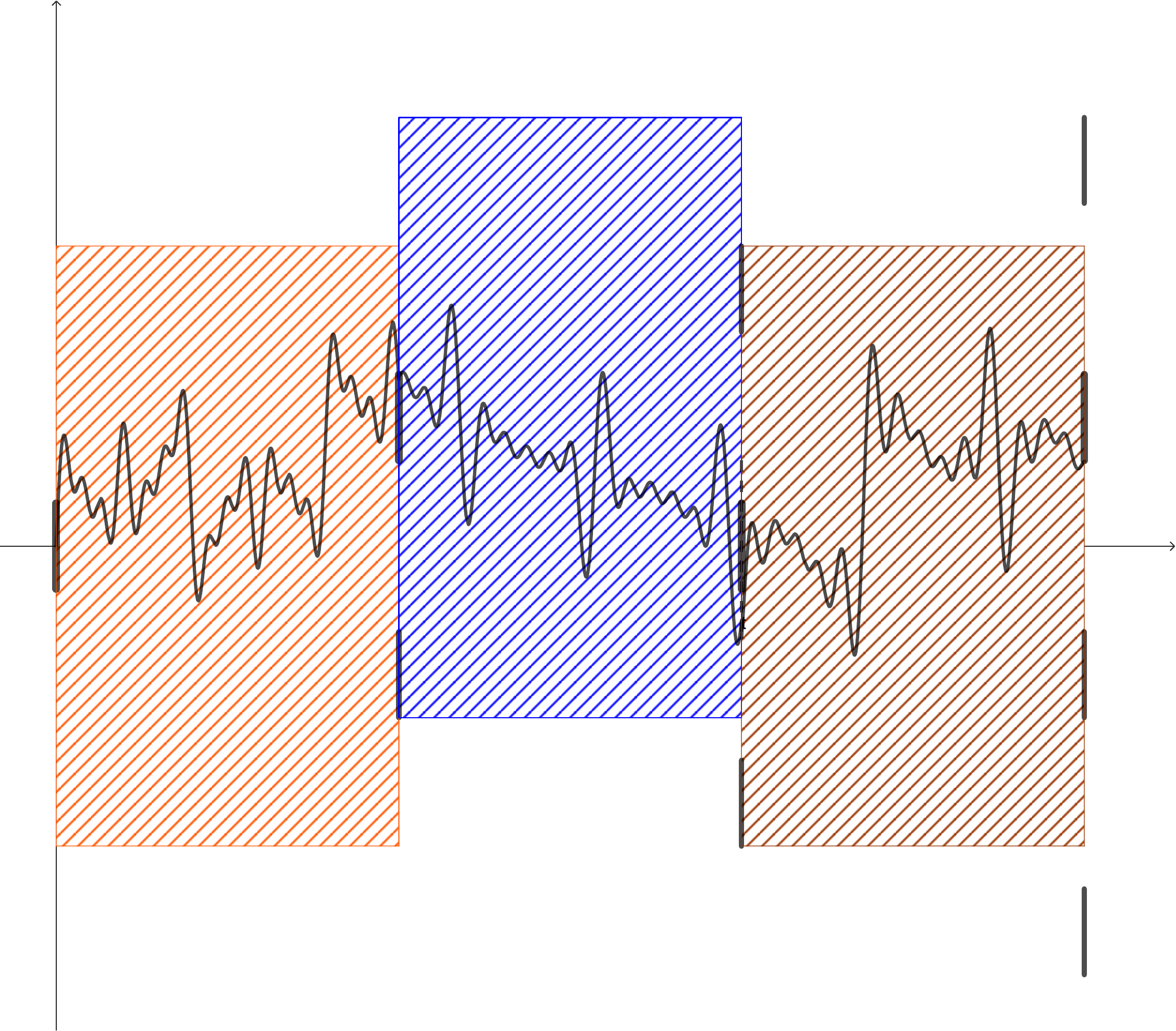}
\caption{Image of SRW trace in tubes: Bold horizontal segments are CG lattice sites. Simple random walk starts from the origin and goes through the site in CG lattice sites without escaping diagonal line areas.}
\end{center}
\end{figure}

Hereafter, we assume  $\beta\in (0,1]$ and we denote by \begin{align*}
n=n_\beta=\lfloor \beta^{-4}\rfloor.
\end{align*}


We obtain a lower bound of the free energy by the following lemma:
\begin{lem}\label{lem:path}

Taking $L>0$ large enough, then for any $\e>0$, there exist  $\delta>0$,  $p>0$, $0<\beta_0\leq 1$, $T_0\geq 1$ such that for any $T\geq T_0$, $\beta\in (0,\beta_0]$ \begin{align*}
Q\left(\begin{array}{l}
\dis \text{There exists a random infinite path in $\mathbb{L}_{CG}$, $\Gamma$, s.t.~}\\
\dis W_{\beta,ITn}\left(\bigcap_{0\leq J\leq I-1}\Omega_{(J,\Gamma_J,\Gamma_{J+1})}^{(T,n,\delta,L)}\right)\geq e^{TI(F_\kZ(\sqrt{2})-\e)} \text{ for all $I\geq 0$}
\end{array}
\right)>p,
\end{align*}
where we define \begin{align*}
W_{\beta,n}(A)=P_S\left[\prod_{k=1}^{n}\zeta_{k,S_{k}}(\beta,\eta): A\right]
\end{align*}
\text{for }$A\in \mathcal{F}_S$, and an infinite path $\Gamma$ is an infinite sequence of sites  $\{(J,\Gamma_J):J\geq 0\}$ with $\Gamma_0=0$, $|\Gamma_{J+1}-\Gamma_J|=1$ $($$J\geq 0$$)$.
\end{lem}

\begin{proof}[Proof of \eqref{eq:lower}] It is trivial that for any infinite path $\Gamma$ \begin{align*}
W_{\beta,ITn}\geq W_{\beta,ITn}\left(\bigcap_{0\leq J\leq I-1}\Omega_{(J,\Gamma_J,\Gamma_{J+1})}^{(T,n,\delta,L)}\right)
\end{align*}
and hence Lemma \ref{lem:path} implies that for any $\e>0$, \begin{align*}
Q\left(\varliminf_{I\to \infty}\frac{1}{IT}\log W_{\beta,ITn}\geq F_\kZ(\sqrt{2})-\e\right)>p
\end{align*}
for large $T\geq 1$ and small $\beta>0$.

Thus, we have \begin{align*}
nF(\beta)\geq F_\kZ(\sqrt{2})-\e
\end{align*}
with positive probability. 
Since the convergence is  $Q$-a.s.~in \eqref{free}, it implies that $\dis \varliminf_{\beta\to 0}\frac{F(\beta)}{\beta^4}\geq F_\kZ(\sqrt{2})-\e$. Thus, \eqref{eq:lower} follows.

\end{proof}
Therefore, it is enough to prove Lemma \ref{lem:path}.


\subsection{Coarse grained oriented site percolation}\label{2.2}

Fix $\e>0$ sufficiently small. To find an infinite path as in Lemma \ref{lem:path}, we will look for ``good" CG-bonds in the following sense.

\begin{definition}
Fix $\e\in (0,1)$.
We say a space-time CG bond $\lan(I,X),(I+1,Y)\ran$ is $\e$-\textit{good} if 
\begin{align*}
\inf_{x\in B_X^{(T,n,\delta)}}\theta_{ITn}\circ W^x_{\beta,Tn}\left(\Omega_{(0,X,Y)}\right)\geq e^{T(F_\kZ(\sqrt{2})-\e)},
\end{align*}
where we define  \begin{align*}
 \theta_n\circ W_{\beta,k}^{x}(A)=P_S^x\left[\prod_{l=1}^{k}\zeta_{l+n,S_l}(\beta,\eta):A\right]
 \end{align*}
 for $k,n\in \N_0$, $x\in \Z$, $A\subset \kF_S$ and infimum is taken over all $x\in B_X^{(T,n,\del)}$ with $x-TIn\in 2\Z$.

Also, we say a space-time CG-bond $\lan(I,X),(I+1,Y)\ran$ is $\e$-\textit{bad} otherwise. 
 \end{definition}

\begin{rem}
Roughly, the bond $\lan (I,X),(I+1,Y)\ran$ is good if the partition functions starting from any space-time sites in $\{ITn\}\times B_{X}$ and restricted to the tube $\mathtt{T}_{(I,X,Y)}$ are sufficiently large.

By the construction and the Markov property  of partition function, we  find that for any infinite path $\Gamma$ and $I\in \N$ \begin{align*}
W_{\beta,TIn}\left(\bigcap_{0\leq J\leq I-1}\Omega_{(J,\Gamma_J,\Gamma_{J+1})}^{(T,n,\delta,L)}\right)\geq \prod_{J=0}^{I-1}\inf_{x\in B_{\Gamma(J)}^{(T,n,\delta)}}\theta_{JTn}\circ W^x_{\beta,Tn}\left(\Omega_{(0,\Gamma_J,\Gamma_{J+1})}^{(T,n,\delta,L)}\right).
\end{align*} 
Thus, Lemma \ref{lem:path} follows when we prove that with positive probability, there exists an infinite path $\Gamma$ such that  each bond $\langle (I,\Gamma_I),(I+1,\Gamma_{I+1})\rangle$ is $\e$-good:

\end{rem}

\begin{lem}\label{lem:OPDPRE}
Take $L>0$ large enough. Then, for any $\e>0$,  there exist  $p_1>0$, $\delta\in (0,1)$, $\beta_0\in (0,1]$, $T_0\geq 1$ such that for any $T\geq T_0$ and $\beta\in (0,\beta_0]$ \begin{align*}
Q\left(\text{There exists an infinite good path $\Gamma$}\right)>p_1.
\end{align*} 
\end{lem}

To obtain an infinite good path, we reduce the problem to the usual oriented percolation. First, we remark that \begin{align*}
&\text{the good-bad configurations of bond }\lan (I,X),(I+1,Y)\ran \text{ and }\lan (J,X'),(J+1,Y')\ran\text{ are idependent}\\
&\text{ if }I\not=J\text{ or }|X-X'|> \frac{2L}{\delta}.
\end{align*}
So it is not the standard oriented percolation, but we use the lemma to compare the good-bad configurations with some supercritical oriented bond percolation.

\begin{lem}\label{lem:stodom1}{\cite[Theorem 1.3]{LSS97}, \cite[Theorem 3.12]{AY15}}
Let $\{X_s\}_{s\in \Z}$ be a collection of $0$-$1$ valued and $k$-dependent random variables, and suppose that there exists a $p\in (0,1)$ such that for each $s\in \Z$\begin{align*}
P\left(X_s=1\right)\geq p.
\end{align*}
Then if \begin{align*}
p>1-\frac{k^k}{(k+1)^{k+1}},
\end{align*}
then $\{X_s\}_{s\in \Z}$ is dominated from below by a product random fields with density $0<\rho(p)<1$. Furthermore, $\rho(p)\to 1$ as $p\to 1$. 
\end{lem}

Let  $U^{(T,n,\delta,L)}_{(I,X,Y)}(\e)=1\{\langle (I,X),(I+1,Y)\rangle \text{ is $\e$-good}\}$ for $\langle (I,X),(I+1,Y)\rangle \in\mathbb{L}_{CG}$. 
Then, $\left\{U^{(T,n,\delta,L)}_{(I,X,Y)}(\e):\langle (I,X),(I+1,Y)\rangle \in \mathbb{L}_{CG}\right\}$ are $0$-$1$ valued and $\frac{4L}{\delta}$-dependent as mentioned above.  
\begin{lem}\label{lem:Berdens}
Take $L>0$ large enough. Then, for any $\e>0$
 there exists   $\delta\in (0,1)$ 
 such that 
 \begin{align*}
\varliminf_{T\to\infty}\varliminf_{\beta\to 0}\inf_{I,X} Q\left(U^{(T,n,\delta,L)}_{(I,X,Y)}(\e )=1\right)=1.
\end{align*}
\end{lem}





\begin{proof}[Proof of Lemma \ref{lem:OPDPRE}]
Lemma \ref{lem:stodom1}  and Lemma \ref{lem:Berdens} yield that $\{U_{(I,X,Y)}^{(T,n,\delta,L)}(\e):\lan (I,X),(I+1,Y)\ran\in\mathbb{L}_{CG}\}$ dominates oriented bond percolation  on $\mathbb{N}_0\times \Z$ with density $0<\rho(p)<1$ from above.   

In particular, we know the critical probability $\overrightarrow{p}$ of oriented bond site percolation is non-trivial.
 Thus, taking $T\geq 1$ large enough and $\beta\in (0,1]$ small enough  such that $Q\left(U_{(I,X,Y)}^{(T,n,\delta,L)}(\e)=1\right)>p$ with $\rho(p)>\overrightarrow{p}$, we have that \begin{align*}
Q\left(\textrm{There exists an infinite nearest neighbor path $\Gamma$}\right)>0.
\end{align*}

\end{proof}

Thus, we complete the proof of Theorem \ref{thm:main} when Lemma \ref{lem:Berdens} is proved.


Hereafter, we omit  $(T,n,\delta,L)$ and $(T,n,\delta)$ if they are clear from the context.

\section{Proof of Lemma \ref{lem:Berdens}}
\subsection{Continuum directed polymers from DPRE}

In this subsection,  we give an idea to prove Lemma \ref{lem:Berdens}. It is clear that the $\{U_{(I,X,Y)}(\e): \langle (I,X),(I+1,Y)\rangle\in \mathbb{L}_{CG}\}$ is identically distributed, so it is enough to focus on $U_{(0,0,1)}(\e)$.

 Recalling the definition of $U$, we need focus on the partition function on paths restricted to the Brownian scaled tube. It is similar to the partition function discussed in  \cite{AKQ14a}, but the weak convergence on path space does not seem to be proved and we will avoid discussing the weak convergence on path space.

The basic idea is the following proposition which is a corollary of \cite[Theorem 2.2, Lemma B.3]{AKQ14a}.

\begin{prop}\label{thmref} Suppose $d=1$. Let $\{\beta=\beta_n:n\geq 1\}$ be an $\mathbb{R}$-valued sequence with $\beta_n=\frac{1}{n^\frac{1}{4}}$. Then,  $\dis \left\{W^{x\sqrt{n}}_{\beta, Tn}\left(\left\{S_{Tn}\in B_{1}^{(T,n,\delta)}\right\}\right):x\in [-1,1]\right\}$ converges weakly (under the topology of the supremum norm)  to a random variable \begin{align}
\int_{\delta T-1}^{\delta T+1}\kZ_{\sqrt{2}}^x(T,y)dy\label{eq:kpzseg}
\end{align}
\end{prop}
\begin{rem}
In \cite[Lemma B.3]{AKQ14a}, $C([\ve,T]\times\R)$-tightness  of the two-parameter field $\left\{\sqrt{n}W_{\beta, sn}\left(\left\{S_{ sn}= x\sqrt{n}\right\}\right):(s,x)\in [\ve,T]\times \R\right\}$ for any $0<\ve<T<\infty$ are proved.
 However, they also discussed  $C(\Delta^2_\ve([0,T])\times\R^2)$-tightness of the four-parameter fields \begin{align*}
\left\{\sqrt{n}\theta_{sn}\circ W^{ x\sqrt{n}}_{\beta,tn}\left(\left\{S_{(t-s)n}= y\sqrt{n}\right\}\right):(s,x),(t,y)\in [0,T]\times \R, t-s\geq \ve\right\}
\end{align*}
 in the second remark in Appendix B in \cite{AKQ14a}, where $\Delta_\ve^2[0,T]=\{(s,t)\in [0,T]^2: t-s\geq \ve\}$. Then, the Riemanian summation \begin{align*}
 W^{x\sqrt{n}}_{\beta, Tn}\left(\left\{S_{Tn}\in B_{1}^{(T,n,\delta)}\right\}\right)=\frac{2}{\sqrt{n}}\sum_{\begin{smallmatrix}y\sqrt{n}\in  B_{1}^{(T,n,\delta)}:\\Tn-y\sqrt{n}\in 2\Z\end{smallmatrix}}\frac{\sqrt{n}}{2}W^{x\sqrt{n}}_{\beta, Tn}\left(\left\{S_{Tn}=y\sqrt{n}\right\}\right)
 \end{align*}
 converges weakly (under the topology of the uniform norm) to \eqref{eq:kpzseg}  by \cite[Theorem 2.2]{AKQ14a} as $\beta\to 0$ (which implies $n\to\infty$).
\end{rem}

We observe that \begin{align*}
&Q\left(\inf_{x\in B_{0}} W^x_{\beta,Tn}\left(\Omega_{(0,0,1)}\right)\geq e^{T(F_\kZ(\sqrt{2})-\e )}\right)\\
&\geq Q\left(\inf_{x\in B_{0}} W^x_{\beta,Tn}\left(\left\{S_{Tn}\in B_1\right\}\right)> 2e^{T(F_\kZ(\sqrt{2})-\e )}\right)\\
&-Q\left(\sup_{x\in B_{0}} W^x_{\beta,Tn}\left(\left\{\{s_i\}_{i=0}^{Tn}:s\not \subset \mathtt{T}_{(0,0)}\right\}\right)> e^{T(F_\kZ(\sqrt{2})-\e )}\right).
\end{align*}
Also, Proposition \ref{thmref} implies that \begin{align*}
\inf_{x\in B_{0}} W^x_{\beta,Tn}\left(\left\{S_{Tn}\in B_1\right\}\right)\Rightarrow \inf_{x\in [-1,1]}\int_{\delta T-1}^{\delta T+1}\kZ_{\sqrt{2}}^x(T,y)dy.
\end{align*}
Therefore, we have that 
\begin{align*}
&\varliminf_{\beta\to0}
Q\left(\inf_{x\in B_{0}^{(T,n,\delta)}} W^x_{\beta,Tn}\left(\Omega_{(0,0,1)}\right)\geq e^{T(F_\kZ(\sqrt{2}-\e ))}\right)\\
&\geq Q\left(\inf_{x\in [-1,1]}\int_{\delta T-1}^{\delta T+1}\kZ_{\sqrt{2}}^x(T,y)dy>2e^{T(F_\kZ(\sqrt{2})-\e)}\right)\\
&-\varlimsup_{\beta\to0}Q\left(\sup_{x\in B_{0}^{(T,n,\delta)}} W^x_{\beta,Tn}\left(\left\{\{s_i\}_{i=0}^{Tn}:s\not \subset \mathtt{T}_{(0,0)}\right\}\right)> e^{T(F_\kZ(\sqrt{2})-\e )}\right).
\end{align*}

The proof is completed when we prove the following two lemmas:
\begin{lem}\label{lem:wout}

Take $L>0$ large enough. Then, for any $\e>0$ and $\delta\in (0,1)$,  
 \begin{align*}
\varlimsup_{T\to\infty}\varlimsup_{\beta\to 0}Q\left(\sup_{x\in B_{0}^{(T,n,\delta)}} W^x_{\beta,Tn}\left(\left\{\{s_i\}_{i=0}^{Tn}:s\not \subset \mathtt{T}_{(0,0)}\right\}\right)> e^{T(F_\kZ(\sqrt{2})-\e )}\right)=0
\end{align*}

\end{lem}

\begin{lem}\label{lem:wexp}
For any $\e>0$,  there exists $\delta \in (0,1)$ such that \begin{align*}
\varliminf_{T\to \infty}Q\left(\inf_{x\in [-1,1]}\int_{\delta T-1}^{\delta T+1}\kZ_{\sqrt{2}}^x(T,y)dy>2e^{T(F_\kZ(\sqrt{2})-\e)}\right)=1
\end{align*}

\end{lem}

\subsection{Proof of Lemma \ref{lem:wout} and Lemma \ref{lem:wexp}}

We will use Garsia-Rodemich-Rumsey's lemma \cite[Lemma A.3.1]{Nua09}.
\begin{lem}\label{Gar}
Let $\phi:[0,\infty)\to [0,\infty)$ and $\Psi:[0,\infty)\to [0,\infty)$ be continuous and stricctly increasing functions satisfying\begin{align*}
\phi(0)=\Psi(0)=0,\ \lim_{t\to \infty}\Psi(t)=\infty.
\end{align*}
Let $f:\mathbb{R}\to \mathbb{R}$ be a continuous function. Provided \begin{align*}
\Gamma=\int_{[-1,1]}\int_{[-1,1]}\Psi\left(\frac{|f(t)-f(s)|}{\phi(|t-s|)}\right)dsdt<\infty,
\end{align*}
then for all $s,t\in [-1,1]$,\begin{align*}
|f(t)-f(s)|\leq 8\int_0^{2|t-s|}\Psi^{-1}\left(\frac{4^{2}\Gamma}{\lambda u^{2}}\right)\phi(du),
\end{align*}
where $\lambda$ is a universal constant.
\end{lem}

Applying Lemma \ref{Gar} with $\Psi(x)=|x|^p$, $\phi(u)=u^q$ for $p\geq 1$, $q>0$ and $pq>2$, we have that \begin{align}
|f(x)-f(y)|\leq \frac{2^{\frac{2}{p}+q+3}}{\lambda^{\frac{1}{p}}\left(q-\frac{2}{p}\right)}|x-y|^{q-\frac{2}{p}}\left(\int_{[-1,1]}\int_{[-1,1]}\left(\frac{|f(t)-f(s)|}{|t-s|^{q}}\right)^pdsdt\right)^{\frac{1}{p}}\label{Garpol}
\end{align}
for $x,y\in [-1,1]$.

Take $f$ as a $C([-1,1])$-valued random variable $X(x)$. Then, the integral in the righthand side of \eqref{Garpol} is independent of variables $x,y$ so we have \begin{align}
E\left[\sup_{|x-y|\leq \delta, x,y\in[-1,1]}|X(x)-X(y)|\right]\leq C_{p,q}\delta^{q-\frac{2}{p}}\left(\int_{-1}^1\int_{-1}^1 \frac{E\left[\left|X(x)-X(y)\right|^p\right]}{|x-y|^{pq}}\right)^\frac{1}{p}\label{eq:GarsiaApp}
\end{align} 
for $p\geq 1$, $q>0$ with $pq>2$.

\begin{proof}[Proof of Lemma \ref{lem:wout}]
We write \begin{align}
E=E(S)=\left\{\{s_i\}_{i=0}^{Tn}:s\not \subset \mathtt{T}_{(0,0)}\right\}.\label{eq:eventE}
\end{align}
for a path $S$.

First, we remark from reflection principle and Chernoff bounds for simple random walk that there exists $I:[0,\infty)\to [0,\infty]$ such that 
\begin{align*}
I(x)\to \infty , x\to \infty
\end{align*}
and 
\begin{align}
Q\left[W_{\beta,Tn}^0 \left(E(S)\right)\right]&=P^0_S\left(E(S)\right)\leq C\exp\left(-I(L)T\right) \label{eq:LDP}
\end{align}
 for constant $C>0$. Hence, taking $L>0$ large enough, we have \begin{align*}
&Q\left( W^0_{\beta,Tn}\left(E(S)\right)> e^{T(F_\kZ(\sqrt{2})-\e )}\right)\leq e^{-I(L)T}
 \end{align*}

To apply \eqref{eq:GarsiaApp} to $X(x)=W_{\beta,Tn}^{\sqrt{n}x}\left(\left\{S_{Tn}\in B_1\right\}\cap \Omega^c\right)$, we will estimate \begin{align*}
Q\left[\left|W_{\beta,Tn}^{\sqrt{n}x}\left(E(S)\right)-W_{\beta,Tn}^{\sqrt{n}y}\left(E(S)\right)\right|^p\right]
\end{align*}
 for $p\geq 1$. To do it, we use the hypercontractivity argument as \cite[Lemma 2.5]{Nak19} developed in \cite{MOO05}. 
 
 
We have \begin{align*}
W_{\beta,Tn}^{\sqrt{n}x}(E)&=P_S^{\sqrt{n}x}\left[\prod_{k=1}^{Tn}\zeta_{k,S_k}:E\right]\\
&=P_S^{\sqrt{n}x}\left[1_{E}\prod_{k=1}^{Tn}(1+(\zeta_{k,S_k}-1))\right]\\
&=\sum_{k=0}^{Tn}\sum_{x_1,\cdots,x_k\in \Z}\sum_{1\leq i_1<\cdots<i_k\leq T_n}\prod_{j=1}^k (\zeta_{i_j,{x_{j}}}-1)P_{S}^{\sqrt{n}x}(S_{i_j}=x_j, j=1,\cdots,k,E).
\end{align*} 
When we set \begin{align*}
\Theta^{(k)}(x)=\begin{cases}
P_S^x(E),\quad &k=0\\
\dis \sum_{x_1,\cdots,x_k\in \Z}\sum_{1\leq i_1<\cdots<i_k\leq Tn}\prod_{j=1}^k (\zeta_{i_j,{x_{j}}}-1)P_{S}^{\sqrt{n}x}(S_{i_j}=x_j, j=1,\cdots,k,E),\quad &k\geq 1,
\end{cases}
\end{align*} 
it is easy to see that \begin{align*}
&Q[\Theta^{(k)}(x)]=0,\quad x\in \Z\\
&Q[\Theta^{(k)}(x)\Theta^{(l)}(y)]=0,\quad k\not=l, x,y\in \Z.
\end{align*}

Then, we have \begin{align*}
&Q\left[|W_{\beta,Tn}^{x}(E)-W_{\beta,Tn}^y(E)|^2\right]\\
&=\sum_{k=0}^n Q\left[\left(\Theta^{(k)}(x)-\Theta^{(k)}(y)\right)^2\right]\\
&=(P_S^x(E)-P_S^y(E))^2\\
&+\sum_{k=1}^{Tn}\sum_{x_1,\cdots,x_k\in \Z}\sum_{1\leq i_1<\cdots<i_k\leq Tn}\\
&\hspace{3em}\prod_{j=1}^k Q\left[(\zeta_{i_j,{x_{j}}}-1)^2\right]\left(P_{S}^{\sqrt{n}x}(S_{i_j}=x_j, j=1,\cdots,k,E)-P_{S}^{\sqrt{n}y}(S_{i_j}=x_j, j=1,\cdots,k,E)\right)^2\\
&=\sum_{k=0}^{Tn} \Lambda^{(k)}(x,y).
\end{align*}
Hypercontractivity developped in \cite[Proposition 3.16 and Proposition 3.12]{MOO05} impies that \begin{align*}
Q\left[|W_{\beta,Tn}^{x}(E)-W_{\beta,Tn}^y(E)|^p\right]\leq \left(\sum_{k=0}^{Tn}\kappa_p^{k}\left(\Lambda^{(k)}(x,y)\right)^{\frac{1}{2}}\right)^p
\end{align*}
where $\dis \kappa_p=2\sqrt{p-1}\sup_{n\geq 1}\frac{Q[\left|\zeta_{1,1}-1\right|^p]^\frac{1}{p}}{Q[\left|\zeta_{1,1}-1\right|^2]^\frac{1}{2}}<\infty $.

From Proposition \ref{prop:i1} and \ref{prop:i3}, there exists $C_3,C_4,C_5,C_6>0$ such that for any $k\geq 1$, $T\geq 0$, $n\geq 1$, $L\geq 2$, and $x,y\in \Z$ \begin{align*}
\Lambda^{(k)}(x,y)\leq C_6\left(C_4\sqrt{T}+\frac{C_5}{LT}\right)|x-y|e^{-I(L-1)T}\frac{(C_3T)^\frac{k-1}{2}}{\Gamma\left(\frac{k-1}{2}+1\right)},
\end{align*}
where  $\Gamma(x)=\int_0^\infty t^{x-1}e^{-t}dt$ is the Gamma funtion. Thus, there exists $C_7>0$ such that 
\begin{align*}
Q\left[|W_{\beta,Tn}^{x}(E)-W_{\beta,Tn}^y(E)|^p\right]&\leq \left( \sum_{k=0}^{Tn}\kappa_p^k\left(\Lambda^{(k)}(x,y)\right)^{\frac{1}{2}}\right)^p\\
&\leq C_{T,p}|x-y|^\frac{p}{2}e^{-\frac{p}{2}I(L-1)T}\left(\sum_{k=0}^{Tn}\frac{1}{2^k}\right)^\frac{p}{2}\left( \sum_{k=0}^{Tn}\frac{(4\kappa_p^2C_3T)^\frac{k-1}{2}}{\Gamma\left(\frac{k-1}{2}+1\right)}\right)^\frac{p}{2},
\end{align*}
where $C_{T,p}$ is a constant which depends on $T$, $p$ (and constants $C_3,\dots,C_6$) and which grows polynomially in $T$. We remark that \begin{align*}
C_7:=\varlimsup_{T\to \infty}\frac{1}{T}\log \sum_{k=0}^{\infty}\frac{(4\kappa_p^2C_3T)^\frac{k-1}{2}}{\Gamma\left(\frac{k-1}{2}+1\right)}<\infty.
\end{align*}

Applying \eqref{eq:GarsiaApp} $p=4$, $q=\frac{5}{8}$, it follows that there exists $C_8>0$ such that \begin{align*}
Q\left[\sup_{|x-y|\leq 1, x,y\in [-1,1]}\left|W_{\beta,Tn}^{\sqrt{n}x}\left(E(S)\right)-W_{\beta,Tn}^{\sqrt{n}y}\left(E(S)\right)\right|\right]\leq C_8(T,p)e^{-\frac{I(L-1)}{2}T}\left(\sum_{k=0}^{\infty}\frac{(4\kappa_p^2C_3T)^\frac{k-1}{2}}{\Gamma\left(\frac{k-1}{2}+1\right)}\right)^\frac{1}{4},
\end{align*}
where $C_{8}(T,p)$ grows at polynomially in $T$. Thus, the statement holds by taking $L>1$ large enough so that \begin{align*}
-\frac{I(L-1)}{2}+C_7<F_{\mathcal{Z}}(\sqrt{2})-\e.
\end{align*}

\end{proof}

\begin{proof}[Proof of Lemma \ref{lem:wexp}]
Translation invariance of white noise and Girsanov transformation of Brownian motion allow us to have \begin{align*}
\left\{\int_{\delta T-1}^{\delta T+1}\mathcal{Z}_{\sqrt{2}}^x(T,y)\dd y:x\in [-1,1]\right\}\stackrel{d}{=}\left\{e^{\delta-\frac{\delta^2T}{2}}\int_{-1}^{1}\mathcal{Z}_{\sqrt{2}}^x(T,y)\dd y:x\in [-1,1]\right\}.
\end{align*}
So, it is enough to show that for some $\delta>0$ \begin{align*}
\varliminf_{T\to \infty}Q\left(\inf_{x\in [-1,1]}\int_{-1}^{1}\kZ_{\sqrt{2}}^x(T,y)dy>2e^{-\delta+\frac{\delta^2T}{ 2}}e^{T(F_\kZ(\sqrt{2})-\e)}\right)=1.
\end{align*}
Also, by taking  $\delta=\sqrt{\e}>0$, it is enough to prove that \begin{align*}
\varliminf_{T\to \infty}Q\left(\inf_{x\in [-1,1]}\int_{-1}^{1}\kZ_{\sqrt{2}}^x(T,y)dy>e^{T(F_\kZ(\sqrt{2})-\frac{\e}{2})}\right)=1
\end{align*} 
We set $\widetilde{\mathcal{Z}}^x(T):=\dis \int_{-1}^1\kZ_{\sqrt{2}}^x(T,y)dy$. Then, modification of the proof in \cite[Lemma 4.4]{Nak19} allows us to see that there exists a $K>0$ such that  for any $\theta \in (-1,1)$, $T>0$
\begin{align*}
Q\left[\exp\left(\theta \left(\log \widetilde{\mathcal{Z}}^x(T)-Q\left[\log \widetilde{\mathcal{Z}}^x(T)\right]\right) \right)\right]\leq \exp\left(-\frac{T\theta^2 K}{1-|\theta|}\right)
\end{align*}

The same argument as in \cite[Corollary 4.5]{Nak19} can be applied to show that \begin{align*}
\frac{1}{T}\log \int_{-1}^1 \mathcal{Z}^0_{\sqrt{2}}(T,y)\dd y \to F_\mathcal{Z}(\sqrt{2}).
\end{align*}
in probability. 

Since \begin{align*}
Q\left(	\inf_{x\in[-1,1]}\widetilde{\mathcal{Z}}^x(T)	\leq e^{T(F_\kZ(\sqrt{2})-\frac{\e}{2})}	\right)&\leq Q\left(\widetilde{\mathcal{Z}}^x(T)\leq 2e^{T(F_\kZ(\sqrt{2})-\frac{\e}{2})} \right)\\
&+Q\left(\sup_{x,y\in[-1,1]}\left|\left(\widetilde{\mathcal{Z}}^x(T)\right)^{-1}-\left(\widetilde{\mathcal{Z}}^y(T)\right)^{-1}\right|\geq \frac{1}{2}e^{-T(F_\kZ(\sqrt{2})+\frac{\e}{2})}\right),
\end{align*}
it is enough to show that \begin{align*}
\varlimsup_{\theta\searrow 0}\varlimsup_{T\to\infty}\frac{1}{\theta T}\log Q\left[\sup_{x,y\in[-1,1]}\left|\left(\widetilde{\mathcal{Z}}^x(T)\right)^{-1}-\left(\widetilde{\mathcal{Z}}^y(T)\right)^{-1}\right|^\theta\right]\leq -F_\mathcal{Z}(\sqrt{2}).
\end{align*}

By H\"older's inequality, we have \begin{align*}
Q\left[\left|\left(\widetilde{\mathcal{Z}}^x(T)\right)^{-1}-\left(\widetilde{\mathcal{Z}}^y(T)\right)^{-1}\right|^\theta \right]\leq Q\left[\left|\widetilde{\mathcal{Z}}^x(T)-\widetilde{\mathcal{Z}}^y(T)\right|^{p\theta }\right]^\frac{1}{p}
Q\left[\left(\widetilde{\mathcal{Z}}^x(T)\right)^{-2q\theta}\right]^\frac{1}{2q}
Q\left[\left(\widetilde{\mathcal{Z}}^y(T)\right)^{-2q\theta}\right]^\frac{1}{2q}
\end{align*}
for any $p,q>1$ with $\frac{1}{p}+\frac{1}{q}=1$.

Applying this inequality to \eqref{eq:GarsiaApp}, we have \begin{align*}
&Q\left[\sup_{x,y\in[-1,1]}\left|\left(\widetilde{\mathcal{Z}}^x(T)\right)^{-1}-\left(\widetilde{\mathcal{Z}}^y(T)\right)^{-1}\right|^\theta\right]\\
&\leq  C\left(\int_{-1}^1\int_{-1}^1 \dd x\dd y \frac{Q\left[\left|\widetilde{\mathcal{Z}}^x(T)-\widetilde{\mathcal{Z}}^y(T)\right|^{p\theta r}\right]^\frac{1}{p}
Q\left[\left(\widetilde{\mathcal{Z}}^x(T)\right)^{-2q\theta r}\right]^\frac{1}{2q}
Q\left[\left(\widetilde{\mathcal{Z}}^y(T)\right)^{-2q\theta r}\right]^\frac{1}{2qr}}{|x-y|^{s}}\right)^\frac{1}{r}
\end{align*}
for $r\geq 1,s>0$ with $rs>2$.

Using the Girasnov transformation again, we can find that \begin{align*}
\sup_{x\in[-1,1]}Q\left[\left(\widetilde{\mathcal{Z}}^x(T)\right)^{-2q\theta r}\right]^\frac{1}{2qr}\leq CQ\left[\left(\widetilde{\mathcal{Z}}^0(T)\right)^{-2q\theta r}\right]^\frac{1}{2qr}.
\end{align*}
and \begin{align*}
\varlimsup_{\theta\searrow 0}\varlimsup_{T\to \infty}\frac{1}{T\theta}\log \sup_{x\in[-1,1]}Q\left[\left(\widetilde{\mathcal{Z}}(T)^x(T)\right)^{-2q\theta r}\right]^\frac{1}{2qr}\leq -2F_{\mathcal{Z}}(\sqrt{2}).
\end{align*}

Also, the same argument as in the proof of \cite[Lemma 4.7]{Nak19} can be applied to obtain 
\begin{align*}
\varlimsup_{\theta\searrow 0}\varlimsup_{T\to \infty}\frac{1}{T\theta}\log \int_{-1}^1\int_{-1}^1 \dd x\dd y \frac{Q\left[\left|\widetilde{\mathcal{Z}}^x(T)-\widetilde{\mathcal{Z}}^y(T)\right|^{p\theta r}\right]^\frac{1}{pr}}{|z-y|^{rs}}\leq F_{\mathcal{Z}}(\sqrt{2})
\end{align*}
for some $r\geq 1$, $s>0$ with $rs>2$.

Thus, the proof is completed.

\end{proof}

\appendix
\section{Local limit theorem}

Now, we  construct  a pair of simple random walks $(S,\widetilde{S})$ on $\Z$ starting from $x$ and $y$ as follows. Let $S$ be a simple random walk starting from $x$ and $\tau_{x,y}=\inf\{n\geq 0:S_n=\frac{x+y}{2}\}$. Then, $\widetilde{S}$ is defined by \begin{align*}
\widetilde{S}_0=y,\quad
\widetilde{S}_n=\begin{cases}
x+y-S_n,\quad &n\leq \tau_{x,y},\\
S_n,\quad &n\geq \tau_{x,y}.
\end{cases}
\end{align*}
Then, $\widetilde{S}$ has the same marginal distribution  as simple random walk starting from $y$. When we denote the law of $(S,\widetilde{S})$ by $P_{S,\widetilde{S}}^{x,y}$, we have that \begin{align*}
|P_S^x(S\in A)-P_S^y(S\in B)|&=\left|P_{S,\widetilde{S}}^{x,y}(S\in A)-P_{S,\widetilde{S}}^{x,y}(\widetilde{S}\in B)\right|\\
&\leq P_{S,\widetilde{S}}^{x,y}(S\in A,\widetilde{S}\in B^c)+P_{S,\widetilde{S}}^{x,y}(S\in A^c,\widetilde{S}\in  B)
\end{align*}
for $A,B$ a measurable set in path space.


Thus, we have 
\begin{align*}
&\left|P_S^{x}(S_{i_j}=x_j,j=1,\cdots,k, E)-P_{S}^y({S}_{j_i}=x_i,i=1,\cdots,k,E)\right|\\
&\leq P_{S,\widetilde{S}}^{x,y}(\{S_{i_j}=x_j,j=1,\cdots,k\}\cap  \{\widetilde{S}_{i_j}=x_j,j=1,\cdots,k\}^c\cap E(S))\\
&+P_{S,\widetilde{S}}^{x,y}(\{S_{i_j}=x_j,j=1,\cdots,k\}\cap E(S)\cap E(\widetilde{S})^c)\\
&+P_{S,\widetilde{S}}^{x,y}(\{S_{i_j}=x_j,j=1,\cdots,k\}^c\cap \{\widetilde{S}_{i_j}=x_j,j=1,\cdots,k\} \cap E(\widetilde{S}))\\
&+P_{S,\widetilde{S}}^{x,y}(\{\widetilde{S}_{i_j}=x_j,j=1,\cdots,k\}\cap E(\widetilde{S})\cap E(S)^c),
\end{align*}
where $E=E(S)$ is the event defined in \eqref{eq:eventE}.

We may assume $0\leq |y|<x$. Then. $E(\widetilde{S})\subset E({S})$ and $E(S)^c\cap E(\widetilde{S})=\emptyset$. Also, this implies that on $E({S})\cap E(\widetilde{S})^c$, $\tau_{x,y}>\tau_{LT\sqrt{n}}$, where $\tau_a=\inf\{n\geq0:|S_n|=a\}$ for $a\in\Z$. Thus, we have \begin{align*}
&\{{S}_{i_j}=x_j,j=1,\cdots,k\}\cap E({S})\cap E(\widetilde{S})^c\\
&\subset \{{S}_{i_j}=x_j,j=1,\cdots,k\}\cap E({S})\cap (\{\tau_{LT\sqrt{n}}< \tau_{x,y} \leq i_1\}\cup \{i_1 \vee \tau_{LT\sqrt{n}}\leq  \tau_{x,y}\})\\
&\subset (\{{S}_{i_j}=x_j,j=1,\cdots,k\}\cap \{\tau_{LT\sqrt{n}}< \tau_{x,y} \leq i_1\})\\
&\quad \cup (\{{S}_{i_j}=x_j,j=1,\cdots,k\}\cap \{\widetilde{S}_{i_1}\not=x_1\}\cap E({S})).
\end{align*}

Hence, we have \begin{align}
&\Lambda^{(k)}(x,y)\notag\\
&\leq 4\frac{2^k}{n^{\frac{k}{2}}}\sum_{x_1,\cdots,x_k\in\Z}\sum_{1\leq j_1<\cdots<j_k\leq Tn}P_{S,\widetilde{S}}^{\sqrt{n}x,\sqrt{n}y}(\{S_{i_j}=x_j,j=1,\cdots,k\}\cap  \{\widetilde{S}_{i_j}=x_j,j=1,\cdots,k\}^c\cap E(S))^2\notag\\
&+ 4\frac{2^k}{n^{\frac{k}{2}}}\sum_{x_1,\cdots,x_k\in\Z}\sum_{1\leq j_1<\cdots<j_k\leq Tn}P_{S,\widetilde{S}}^{\sqrt{n}x,\sqrt{n}y}(\{S_{i_j}=x_j,j=1,\cdots,k\}^c\cap  \{\widetilde{S}_{i_j}=x_j,j=1,\cdots,k\}\cap E(\widetilde{S}))^2\notag\\
&+4\frac{2^k}{n^{\frac{k}{2}}}\sum_{x_1,\cdots,x_k\in\Z}\sum_{1\leq j_1<\cdots<j_k\leq Tn}P_{S,\widetilde{S}}^{\sqrt{n}x,\sqrt{n}y}(\{{S}_{i_j}=x_j,j=1,\cdots,k\}\cap \{\tau_{LT\sqrt{n}}< \tau_{x,y} \leq i_1\})^2\notag\\
&+4\frac{2^k}{n^{\frac{k}{2}}}\sum_{x_1,\cdots,x_k\in\Z}\sum_{1\leq j_1<\cdots<j_k\leq Tn}P_{S,\widetilde{S}}^{\sqrt{n}x,\sqrt{n}y}(\{{S}_{i_j}=x_j,j=1,\cdots,k\}\cap \{\widetilde{S}_{i_1}\not=x_1\}\cap E({S}))^2\notag\\
&\leq 8\frac{2^k}{n^{\frac{k}{2}}}\sum_{x_1,\cdots,x_k\in\Z}\sum_{1\leq j_1<\cdots<j_k\leq Tn}\notag\\
&\hspace{6em}  P_{S,\widetilde{S}}^{\sqrt{n}x,\sqrt{n}y}(\{S_{i_1}=x_j,j=1,\cdots,k\}\cap  \{\widetilde{S}_{i_1}\not=x_1\})P_S^{\sqrt{n}x}(\{S_{i_j}=x_j,j=1,\cdots,k\}\cap E(S))\notag\\
&+ 4\frac{2^k}{n^{\frac{k}{2}}}\sum_{x_1,\cdots,x_k\in\Z}\sum_{1\leq j_1<\cdots<j_k\leq Tn}\notag\\
&\hspace{6em}  P_{S,\widetilde{S}}^{\sqrt{n}x,\sqrt{n}y}(\{S_{i_1}\not=x_1\}\cap  \{\widetilde{S}_{i_j}=x_j,j=1,\cdots,k\})P_{\widetilde{S}}^{\sqrt{n}y}(\{\widetilde{S}_{i_j}=x_j,j=1,\cdots,k\}\cap E(\widetilde{S}))\notag\\
&+4\frac{2^k}{n^{\frac{k}{2}}}\sum_{x_1,\cdots,x_k\in\Z}\sum_{1\leq j_1<\cdots<j_k\leq Tn}P_{S,\widetilde{S}}^{\sqrt{n}x,\sqrt{n}y}(\{{S}_{i_j}=x_j,j=1,\cdots,k\}\cap \{\tau_{LT\sqrt{n}}< \tau_{\sqrt{n}x,\sqrt{n}y} \leq i_1\})^2\notag\\
&=:8I_1+4I_2+4I_3,
 \label{eq:Lambdadiff}
\end{align}
where we have used $Q\left[(\zeta_{1,1}-1)^2\right]\leq \frac{2}{n^{\frac{1}{2}}}$ for large $n$.

From symmetry, we may focus on only $I_1$ and $I_3$.

\begin{prop}\label{prop:i1}
There exists $C_3>0$ and $C_4>0$ such that for any $k\geq 1$, $T\geq 0$, $n\geq 1$, $L\geq 2$, and $x,y\in \Z$ \begin{align*}
I_1\leq Ce^{I(L-1)}C_4|x-y|\sqrt{T}\frac{\left(C_3T\right)^{\frac{k-1}{2}}}{\Gamma\left(\frac{k-1}{2}+1\right)},
\end{align*}
where $I$ is the function defined in \eqref{eq:LDP}.
\end{prop}

\begin{proof}
It is easy to see from local limit theorem and \eqref{eq:LDP} that \begin{align*}
I_1&\leq \frac{2^k}{n^{\frac{k}{2}}}\sum_{1\leq i_1<\cdots<i_k\leq Tn}
\sup_{x_1,\dots,x_k}P_{S,\widetilde{S}}^{\sqrt{n}x,\sqrt{n}y}(\{S_{i_1}=x_j,j=1,\cdots,k\}\cap  \{\widetilde{S}_{i_1}\not=x_1\})
P_S^{\sqrt{n}x}( E(S))\\
&\leq Ce^{-I(L-1)T}\frac{2^k}{n^{\frac{k}{2}}}\sum_{1\leq i_1<\cdots<i_k\leq Tn}\sup_{x_1,\dots,x_k}P_{S,\widetilde{S}}^{\sqrt{n}x,\sqrt{n}y}(\{S_{i_1}=x_j,j=1,\cdots,k\}\cap  \{\widetilde{S}_{i_1}\not=x_1\})\\
&\leq Ce^{-I(L-1)T}\frac{2^k}{n^{\frac{k}{2}}}\sum_{ 1\leq i_1<\dots<i_k\leq Tn}P_{S,\widetilde{S}}^{x,y}(\{\tau_{\sqrt{n}x,\sqrt{n}y}\geq i_1\})\prod_{j= 2}^k \frac{c}{\sqrt{i_j-i_{j-1}}}.
\end{align*}
Enlarging the range of summations on $i_1,\dots,i_k$ to $1\leq i_1\leq Tn$, and $i_1<i_2<\dots<i_k\leq Tn+i_1$, we have \begin{align*}
I_1&\leq Ce^{-I(L-1)T}\frac{2^k}{n^{\frac{k}{2}}}\sum_{ 1\leq i_1<\dots<i_k\leq Tn}P_{S,\widetilde{S}}^{\sqrt{n}x,\sqrt{n}y}(\tau_{\sqrt{n}x,\sqrt{n}y}\geq i_1)\prod_{j= 2}^k \frac{c}{\sqrt{i_j-i_{j-1}}}
\end{align*}
and it follows from \cite[Section 3 and Lemma A.1]{AKQ14a} that there exists $C_3>0$ such that  \begin{align*}
I_1&\leq Ce^{-I(L-1)T}\frac{2}{n^{\frac{1}{2}}}\sum_{1\leq i\leq Tn}P_{S,\widetilde{S}}^{\sqrt{n}x,\sqrt{n}y}(\tau_{\sqrt{n}x,\sqrt{n}y}\geq i)\frac{\left(C_3T\right)^{\frac{k-1}{2}}}{\Gamma\left(\frac{k-1}{2}+1\right)}.
\end{align*}
Also, the reflection principle of simple random walk and the local limit theorem yield that \begin{align*}
P_{S,\widetilde{S}}^{\sqrt{n}x,\sqrt{n}y}(\tau_{\sqrt{n}x,\sqrt{n}y}\geq i)=P_S^{0}\left(-\frac{x-y}{2}\sqrt{n}<S_i\leq \frac{x-y}{2}\sqrt{n}\right)\leq \frac{C_4|x-y|\sqrt{n}}{\sqrt{i}}\wedge 1.
\end{align*}
By the approximation of the Riemmanian summation, we find that \begin{align*}
I_1\leq Ce^{-I(L-1)T}C_4|x-y|\sqrt{T}\frac{\left(C_3T\right)^{\frac{k-1}{2}}}{\Gamma\left(\frac{k-1}{2}+1\right)}.
\end{align*}
\end{proof}

\begin{prop}\label{prop:i3}
There exists $C_5>0$ such that for any $k\geq 1$, $T\geq 0$, $n\geq 1$, $L\geq 2$, and $x,y\in \Z$ \begin{align*}
I_1\leq Ce^{-I(L-1)T}\frac{C_5|x-y|}{LT}\frac{(C_3T)^\frac{k-1}{2}}{\Gamma\left(\frac{k-1}{2}+1\right)},
\end{align*}
where $C_3>0$ is given in Proposition \ref{prop:i1}.
\end{prop}
\begin{proof}
It is clear from the Markov property and the local limit theorem that 
\begin{align*}
&\sum_{x_1\in \Z}\sup_{x_2,\dots,x_k}P_{S,\widetilde{S}}^{\sqrt{n}x,\sqrt{n}y}(\{{S}_{i_j}=x_j,j=1,\cdots,k\}\cap \{\tau_{LT\sqrt{n}}< \tau_{\sqrt{n}x,\sqrt{n}y} \leq i_1\})\\
&\leq P_{S}^{\sqrt{n}x}(\tau_{LT\sqrt{n}}<\tau_{\frac{x+y}{2}\sqrt{n}})P_S^{LT\sqrt{n}}(\tau_{\frac{x+y}{2}\sqrt{n}}\leq Tn)\prod_{j=2}^k \frac{c}{\sqrt{i_j-i_{j-1}}}\\
&\leq \frac{\frac{x-y}{2}}{LT-\frac{x-y}{2}}e^{-I(L-1)T}\prod_{j=2}^k \frac{c}{\sqrt{i_j-i_{j-1}}}
\end{align*}
and 
\begin{align*}
&\sup_{x_1}\sum_{x_2,\dots,x_k}P_{S,\widetilde{S}}^{\sqrt{n}x,\sqrt{n}y}(\{{S}_{i_j}=x_j,j=1,\cdots,k\}\cap \{\tau_{LT\sqrt{n}}< \tau_{\sqrt{n}x,\sqrt{n}y} \leq i_1\})\\
&\leq \sup_{x_1}P_{S}^{\sqrt{n}x}(S_{i_1}=x_1)\leq \frac{c}{\sqrt{ i_1}}.
\end{align*}

Thus, it follorws from \cite[Section 3 and Lemma A.1]{AKQ14a}\begin{align*}
I_3&\leq \frac{2^k}{n^{\frac{k}{2}}}\sum_{1\leq i_1<\dots<i_k\leq Tn}C\frac{|x-y|}{LT}e^{-I(L-1)T}\frac{c}{\sqrt{i_1}}\prod_{j=2}^k\frac{c}{\sqrt{i_j-i_{j-1}}}\\
&\leq Ce^{-I(L-1)T}\frac{C_5|x-y|}{LT}\frac{(C_3T)^\frac{k-1}{2}}{\Gamma\left(\frac{k-1}{2}+1\right)}.
\end{align*}
\end{proof}





\end{document}